\input amssym.def
\input amssym.tex
\font\titre=cmb10 at 14pt

\def\pointir{\unskip . --- \ignorespaces}


\def\Bigbreak{\vskip-\lastskip\bigbreak}
\def\Medbreak{\vskip-\lastskip\medbreak}


\def\ctexte#1\endctexte{%
  \hbox{$\vcenter{\halign{\hfill##\hfill\crcr#1\crcr}}$}}


\long\def\ctitre#1\endctitre{%
    \ifdim\lastskip<24pt\vskip-\lastskip\bigbreak\bigbreak\fi
  		\vbox{\parindent=0pt\leftskip=0pt plus 1fill
          \rightskip=\leftskip
          \parfillskip=0pt\bf#1\par}
    \bigskip\nobreak}

\long\def\section#1\endsection{%
\vskip 0pt plus 3\normalbaselineskip
\penalty-250
\vskip 0pt plus -3\normalbaselineskip
\Bigbreak
\message{[section \string: #1]}{\bf#1\unskip}\pointir}

\long\def\sectiona#1\endsection{%
\vskip 0pt plus 3\normalbaselineskip
\penalty-250
\vskip 0pt plus -3\normalbaselineskip
\Bigbreak
\message{[sectiona \string: #1]}%
{\bf#1}\medskip\nobreak}

\long\def\subsection#1\endsubsection{%
\Medbreak
{\it#1\unskip}\pointir}

\long\def\subsectiona#1\endsubsection{%
\Medbreak
{\it#1}\par\nobreak}

\def\rem#1\endrem{%
\Medbreak
{\it#1\unskip} : }

\def\remp#1\endrem{%
\Medbreak
{\pc #1\unskip}\pointir}

\def\rema#1\endrem{%
\Medbreak
{\it #1}\par\nobreak}

\def\newparwithcolon#1\endnewparwithcolon{
\Medbreak
{#1\unskip} : }

\def\newparwithpointir#1\endnewparwithpointir{
\Medbreak
{#1\unskip}\pointir}

\def\newpara#1\endnewpar{
\Medbreak
{#1\unskip}\smallskip\nobreak}


\long\def\th#1 #2\enonce#3\endth{%
   \Medbreak
   {\pc#1} {#2\unskip}\pointir{\it #3}\medskip}

\long\def\tha#1 #2\enonce#3\endth{%
   \Medbreak
   {\pc#1} {#2\unskip}\par\nobreak{\it #3}\medskip}

\def\Grille{\setbox13=\vbox to 5\unitlength{\hrule width 109mm\vfill} 
\setbox13=\vbox to 65mm{\offinterlineskip\leaders\copy13\vfill\kern-1pt\hrule} 
\setbox14=\hbox to 5\unitlength{\vrule height 65mm\hfill} 
\setbox14=\hbox to 109mm{\leaders\copy14\hfill\kern-2mm\vrule height 65mm}
\ht14=0pt\dp14=0pt\wd14=0pt \setbox13=\vbox to
0pt{\vss\box13\offinterlineskip\box14} \wd13=0pt\box13}


\def\fleche(#1,#2)\dir(#3,#4)\long#5{%
\noalign{\leftput(#1,#2){\vector(#3,#4){#5}}}}

\def\ligne(#1,#2)\dir(#3,#4)\long#5{%
\noalign{\leftput(#1,#2){\lline(#3,#4){#5}}}}

\def\put(#1,#2)#3{\noalign{\setbox1=\hbox{%
    \kern #1\unitlength
    \raise #2\unitlength\hbox{$#3$}}%
    \ht1=0pt \wd1=0pt \dp1=0pt\box1}}


\def\diagram#1{\def\normalbaselines{\baselineskip=0pt\lineskip=5pt}
\matrix{#1}}

\def\hfl#1#2#3{\smash{\mathop{\hbox to#3{\rightarrowfill}}\limits
^{\scriptstyle#1}_{\scriptstyle#2}}}

\def\gfl#1#2#3{\smash{\mathop{\hbox to#3{\leftarrowfill}}\limits
^{\scriptstyle#1}_{\scriptstyle#2}}}


 \message{`lline' & `vector' macros from LaTeX}
 \catcode`@=11
\def\{{\relax\ifmmode\lbrace\else$\lbrace$\fi}
\def\}{\relax\ifmmode\rbrace\else$\rbrace$\fi}
\def\newcount{\alloc@0\count\countdef\insc@unt}
\def\newdimen{\alloc@1\dimen\dimendef\insc@unt}
\def\newwrite{\alloc@7\write\chardef\sixt@@n}

\newwrite\@unused
\newcount\@tempcnta
\newcount\@tempcntb
\newdimen\@tempdima
\newdimen\@tempdimb
\newbox\@tempboxa

\def\@spaces{\space\space\space\space}
\def\@whilenoop#1{}
\def\@whiledim#1\do #2{\ifdim #1\relax#2\@iwhiledim{#1\relax#2}\fi}
\def\@iwhiledim#1{\ifdim #1\let\@nextwhile=\@iwhiledim
        \else\let\@nextwhile=\@whilenoop\fi\@nextwhile{#1}}
\def\@badlinearg{\@latexerr{Bad \string\line\space or \string\vector
   \space argument}}
\def\@latexerr#1#2{\begingroup
\edef\@tempc{#2}\expandafter\errhelp\expandafter{\@tempc}%
\def\@eha{Your command was ignored.
^^JType \space I <command> <return> \space to replace it
  with another command,^^Jor \space <return> \space to continue without
it.} 
\def\@ehb{You've lost some text. \space \@ehc}
\def\@ehc{Try typing \space <return>
  \space to proceed.^^JIf that doesn't work, type \space X <return> \space to
  quit.}
\def\@ehd{You're in trouble here.  \space\@ehc}

\typeout{LaTeX error. \space See LaTeX manual for explanation.^^J
 \space\@spaces\@spaces\@spaces Type \space H <return> \space for
 immediate help.}\errmessage{#1}\endgroup}
\def\typeout#1{{\let\protect\string\immediate\write\@unused{#1}}}

\font\tenln    = line10
\font\tenlnw   = linew10

\newdimen\@wholewidth
\newdimen\@halfwidth
\newdimen\unitlength 

\unitlength =1pt


\def\thinlines{\let\@linefnt\tenln \let\@circlefnt\tencirc
  \@wholewidth\fontdimen8\tenln \@halfwidth .5\@wholewidth}
\def\thicklines{\let\@linefnt\tenlnw \let\@circlefnt\tencircw
  \@wholewidth\fontdimen8\tenlnw \@halfwidth .5\@wholewidth}

\def\linethickness#1{\@wholewidth #1\relax \@halfwidth .5\@wholewidth}

\newif\if@negarg

\def\lline(#1,#2)#3{\@xarg #1\relax \@yarg #2\relax
\@linelen=#3\unitlength
\ifnum\@xarg =0 \@vline
  \else \ifnum\@yarg =0 \@hline \else \@sline\fi
\fi}

\def\@sline{\ifnum\@xarg< 0 \@negargtrue \@xarg -\@xarg \@yyarg -\@yarg
  \else \@negargfalse \@yyarg \@yarg \fi
\ifnum \@yyarg >0 \@tempcnta\@yyarg \else \@tempcnta -\@yyarg \fi
\ifnum\@tempcnta>6 \@badlinearg\@tempcnta0 \fi
\setbox\@linechar\hbox{\@linefnt\@getlinechar(\@xarg,\@yyarg)}%
\ifnum \@yarg >0 \let\@upordown\raise \@clnht\z@
   \else\let\@upordown\lower \@clnht \ht\@linechar\fi
\@clnwd=\wd\@linechar
\if@negarg \hskip -\wd\@linechar \def\@tempa{\hskip -2\wd\@linechar}\else
     \let\@tempa\relax \fi
\@whiledim \@clnwd <\@linelen \do
  {\@upordown\@clnht\copy\@linechar
   \@tempa
   \advance\@clnht \ht\@linechar
   \advance\@clnwd \wd\@linechar}%
\advance\@clnht -\ht\@linechar
\advance\@clnwd -\wd\@linechar
\@tempdima\@linelen\advance\@tempdima -\@clnwd
\@tempdimb\@tempdima\advance\@tempdimb -\wd\@linechar
\if@negarg \hskip -\@tempdimb \else \hskip \@tempdimb \fi
\multiply\@tempdima \@m
\@tempcnta \@tempdima \@tempdima \wd\@linechar \divide\@tempcnta \@tempdima
\@tempdima \ht\@linechar \multiply\@tempdima \@tempcnta
\divide\@tempdima \@m
\advance\@clnht \@tempdima
\ifdim \@linelen <\wd\@linechar
   \hskip \wd\@linechar
  \else\@upordown\@clnht\copy\@linechar\fi}

\def\@hline{\ifnum \@xarg <0 \hskip -\@linelen \fi
\vrule height \@halfwidth depth \@halfwidth width \@linelen
\ifnum \@xarg <0 \hskip -\@linelen \fi}

\def\@getlinechar(#1,#2){\@tempcnta#1\relax\multiply\@tempcnta 8
\advance\@tempcnta -9 \ifnum #2>0 \advance\@tempcnta #2\relax\else
\advance\@tempcnta -#2\relax\advance\@tempcnta 64 \fi
\char\@tempcnta}

\def\vector(#1,#2)#3{\@xarg #1\relax \@yarg #2\relax
\@linelen=#3\unitlength
\ifnum\@xarg =0 \@vvector
  \else \ifnum\@yarg =0 \@hvector \else \@svector\fi
\fi}

\def\@hvector{\@hline\hbox to 0pt{\@linefnt
\ifnum \@xarg <0 \@getlarrow(1,0)\hss\else
    \hss\@getrarrow(1,0)\fi}}

\def\@vvector{\ifnum \@yarg <0 \@downvector \else \@upvector \fi}

\def\@svector{\@sline
\@tempcnta\@yarg \ifnum\@tempcnta <0 \@tempcnta=-\@tempcnta\fi
\ifnum\@tempcnta <5
  \hskip -\wd\@linechar
  \@upordown\@clnht \hbox{\@linefnt  \if@negarg
  \@getlarrow(\@xarg,\@yyarg) \else \@getrarrow(\@xarg,\@yyarg) \fi}%
\else\@badlinearg\fi}

\def\@getlarrow(#1,#2){\ifnum #2 =\z@ \@tempcnta='33\else
\@tempcnta=#1\relax\multiply\@tempcnta \sixt@@n \advance\@tempcnta
-9 \@tempcntb=#2\relax\multiply\@tempcntb \tw@
\ifnum \@tempcntb >0 \advance\@tempcnta \@tempcntb\relax
\else\advance\@tempcnta -\@tempcntb\advance\@tempcnta 64
\fi\fi\char\@tempcnta}

\def\@getrarrow(#1,#2){\@tempcntb=#2\relax
\ifnum\@tempcntb < 0 \@tempcntb=-\@tempcntb\relax\fi
\ifcase \@tempcntb\relax \@tempcnta='55 \or
\ifnum #1<3 \@tempcnta=#1\relax\multiply\@tempcnta
24 \advance\@tempcnta -6 \else \ifnum #1=3 \@tempcnta=49
\else\@tempcnta=58 \fi\fi\or
\ifnum #1<3 \@tempcnta=#1\relax\multiply\@tempcnta
24 \advance\@tempcnta -3 \else \@tempcnta=51\fi\or
\@tempcnta=#1\relax\multiply\@tempcnta
\sixt@@n \advance\@tempcnta -\tw@ \else
\@tempcnta=#1\relax\multiply\@tempcnta
\sixt@@n \advance\@tempcnta 7 \fi\ifnum #2<0 \advance\@tempcnta 64 \fi
\char\@tempcnta}

\def\@vline{\ifnum \@yarg <0 \@downline \else \@upline\fi}

\def\@upline{\hbox to \z@{\hskip -\@halfwidth \vrule
  width \@wholewidth height \@linelen depth \z@\hss}}

\def\@downline{\hbox to \z@{\hskip -\@halfwidth \vrule
  width \@wholewidth height \z@ depth \@linelen \hss}}

\def\@upvector{\@upline\setbox\@tempboxa\hbox{\@linefnt\char'66}\raise
     \@linelen \hbox to\z@{\lower \ht\@tempboxa\box\@tempboxa\hss}}

\def\@downvector{\@downline\lower \@linelen
      \hbox to \z@{\@linefnt\char'77\hss}}

\thinlines

\newcount\@xarg
\newcount\@yarg
\newcount\@yyarg
\newcount\@multicnt
\newdimen\@xdim
\newdimen\@ydim
\newbox\@linechar
\newdimen\@linelen
\newdimen\@clnwd
\newdimen\@clnht
\newdimen\@dashdim
\newbox\@dashbox
\newcount\@dashcnt
 \catcode`@=12


\newbox\tbox
\newbox\tboxa

\def\leftzer#1{\setbox\tbox=\hbox to 0pt{#1\hss}%
     \ht\tbox=0pt \dp\tbox=0pt \box\tbox}

\def\rightzer#1{\setbox\tbox=\hbox to 0pt{\hss #1}%
     \ht\tbox=0pt \dp\tbox=0pt \box\tbox}

\def\centerzer#1{\setbox\tbox=\hbox to 0pt{\hss #1\hss}%
     \ht\tbox=0pt \dp\tbox=0pt \box\tbox}

%
\def\image(#1,#2)#3{\vbox to #1{\offinterlineskip
    \vss #3 \vskip #2}}


\def\leftput(#1,#2)#3{\setbox\tboxa=\hbox{%
    \kern #1\unitlength
    \raise #2\unitlength\hbox{\leftzer{#3}}}%
    \ht\tboxa=0pt \wd\tboxa=0pt \dp\tboxa=0pt\box\tboxa}

\def\rightput(#1,#2)#3{\setbox\tboxa=\hbox{%
    \kern #1\unitlength
    \raise #2\unitlength\hbox{\rightzer{#3}}}%
    \ht\tboxa=0pt \wd\tboxa=0pt \dp\tboxa=0pt\box\tboxa}

\def\centerput(#1,#2)#3{\setbox\tboxa=\hbox{%
    \kern #1\unitlength
    \raise #2\unitlength\hbox{\centerzer{#3}}}%
    \ht\tboxa=0pt \wd\tboxa=0pt \dp\tboxa=0pt\box\tboxa}

\unitlength=1mm
\let\inj=\hookrightarrow
\let\lgr=\longrightarrow
\def\pc#1#2{#1{\sevenrm #2}}
\magnification=\magstep1

\centerline{\titre A toro\"\i dal resolution for the bad reduction} 
\centerline{\titre of  some  Shimura
varieties} 

\vskip 5mm\centerline{\bf by Alain Genestier}
\vskip 2cm

{\bf 1\pointir Statement of the main result}\vskip 3mm
Let $S_K(G,\mu)$ be a Shimura variety. We assume that $G$ is either the group of unitary
similitudes for a hermitian space defined over an imaginary quadratic extension $E$ of $\Bbb Q$
(unitary case) or a group of symplectic similitudes over $\Bbb Q$ (symplectic case).
Let $p$ be a prime in ${\Bbb Z}$. We assume that the group $K_p$ defining the $p$-part of the
level structures is a parahoric subgroup of $G({\Bbb Q}_p)$. In the symplectic case, we assume
that $K_p$ is contained in the Siegel parahoric.
In the unitary case,  we also assume that the place $p$ is split in $E$: $p=\goth P_0 \goth
P_1$.

Following Rapoport and Zink ([RZ]), under these conditions  the Shimura variety $S_K(G,\mu)$ has a
natural modular model
$S_K(G,\mu)_{{\Bbb Z}_{(p)}}$ over ${\rm Spec\,}{\Bbb Z}_{(p)}$. In fact in the unitary case
this model comes from a scheme over the semi-local ring ${\cal O}_{E,(p)}$ by restriction of
scalars.

In some  cases the singularities of this modular  model  have been studied.  

For modular curves, the modular model is the one constructed  by Deligne and\break
Rapoport ([DR]).
The unitary Shimura varieties defined by a hermitian form of signature $(1,n-1)$
in the archimedean place have been studied by M.  Rapoport ([R]), who  proved that
$S_K(G,\mu)_{{\cal O}_{E,(\goth P_i)}}\ (i=0,1)$ is semi-stable over ${\rm Spec\,}{\cal
O}_{E,(\goth P_i)}$.
 The Siegel modular variety with level structures associated to the Siegel parahoric have been
studied by Cha\"\i\ and Normann ([CN]). They proved that the singularities are Cohen-Macaulay,
and Faltings ([F 1]) even proved that they are rational. 
Deligne and Pappas ([DP]) also have a result for  Hilbert-Blumenthal varieties, when
$p$ divides the discriminant of the defining field (note that these varieties do not satisfy the
requirements of the present work).

However, except in the first two cases, the structure morphism is far from beeing semi-stable
and  we are interested in constructing resolutions of
$S_K(G,\mu)_{{\Bbb Z}_{(p)}}$ with  a ``reasonable structure morphism". By this,
we mean a scheme (or eventually an algebraic space)
$\widetilde S_K(G,\mu)_{{\Bbb Z}_{(p)}}$ over $S_K(G,\mu)_{{\Bbb Z}_{(p)}}$ such that 
the morphism
$\widetilde S_K(G,\mu)_{{\Bbb Z}_{(p)}}\lgr S_K(G,\mu)_{{\Bbb Z}_{(p)}}$ is proper and becomes an
isomorphism over $\Bbb Q$
and the morphism $\widetilde S_K(G,\mu)_{{\Bbb Z}_{(p)}}\lgr{\rm Spec\,}{\Bbb Z}_{(p)}$ 
has well understood singularities  (e.g. is semi-stable, or log-smooth for certain given
log-structures).

Results of this type have already been obtained by several authors:

\item {--} the Siegel modular variety of genus 2 has been dealt with by A. J. de Jong, who
showed ([dJ 1]) that a certain blowing-up of  the modular model is semi-stable over ${\Bbb
Z}_{(p)}$

\item {--} using De Concini-Procesi's compactifications of symmetric spaces (and proving that
they can be performed over $\Bbb Z$), Faltings ([F 1]) constructs a semi-stable resolution  of
the Siegel modular variety with level structures defined by the Siegel parahoric. His methods
also work for an analogous parahoric in the unitary case

\item {--} in a recent preprint ([F 2])  which was a source of inspiration for the present
work,  Faltings resolves the unitary Shimura varieties whose defining hermitian form has
signature
$(2,n-2)$ or
$(3,n-3)$ and the Siegel modular  variety of genus 3

\item {--} a semi-stable resolution of the Siegel modular variety of genus 3 with Iwahori
level-structures is also constructed in [G].

In the present work we prove the following one.
\th T{HEOREM} A
\enonce Under the above assumptions for the groups $G$ and $K_p$, $S_K(G,\mu)_{{\Bbb
Z}_{(p)}}$ has a canonical log-smooth resolution.
\endth

{\it Remark}: In fact the resolution $\widetilde S_K(G,\mu)_{{\Bbb
Z}_{(p)}}$ that we obtain satisfies a stronger
condition (not appealing to logarithmic geometry: {\it cf.} 2.2.6) than log--smoothness over
${\Bbb Z}_{(p)}$. Using this stronger statement and ([KKMS], III. 4), one can see that
$S_K(G,\mu)_{{\Bbb
Z}_{(p)}}$ has a (non canonical) semi--stable  resolution over an extension ${\Bbb
Z}_{(p)}[p^{1/\nu}]$ of ${\Bbb
Z}_{(p)}$.

 To prove theorem A, it will be sufficient to prove an
analogous statement for the  local model. Following Cha\"\i,
Norman, Rapoport and Zink ([CN], [R], [RZ]), $S_K(G,\mu)_{{\Bbb Z}_{(p)}}$ has a local model
${\Bbb M}_p$ ; this is a projective scheme over ${\Bbb Z}_{(p)}$ such that for every point $s$ in
$S_K(G,\mu)_{{\Bbb Z}_{(p)}}$, there exists a (non unique) point $m$ in ${\Bbb M}_p$ such that $s$
and $m$ have isomorphic \'etale neighbourhoods. This scheme is defined in terms of linear
algebra and so is easier to deal with than ${\cal S}(g,\mu)_{{\Bbb Z}_p}$. The couple
$(G,\mu)$  enters in its definition  only through   $G_{{\Bbb Q}_p}$ (i.e.
$({\rm GL}_n\times {\Bbb G}_m)_{{\Bbb Q}_p}
$ or ${\rm Gsp}_{2g,{\Bbb Q}_p}$) and (in the unitary case) the signature $(r, n-r)$ of the
hermitian form  ; the parahoric
$K_p$ enters in its definition  only through the partition of $n$ or $g$
defining it (recall that in the symplectic case we assumed that $K_p$ is contained in the
Siegel parahoric). 

The local model
${\Bbb M}_p$ is endowed with an action of a certain affine smooth ${\Bbb Z}_{(p)}$-group-scheme
$\cal K$ associated to the combinatorial data defining the parahoric $K_p={\cal K}({\Bbb
Z}_p)$ (see [dJ 2], [RZ]), and one has a diagram
$$\diagram{
& {\cal T}
\cr\noalign{\vskip 8mm}
{\cal S}_K(G,\mu)_{{\Bbb Z}_{(p)}}
&\kern 10mm& {\Bbb M}_p
\cr
\fleche(24,14)\dir(-1,-1)\long{6}
\fleche(30,14)\dir(1,-1)\long{6}
}$$
where the morphism ${\cal T}\lgr {\cal S}_K(G,\mu)_{{\Bbb Z}_{(p)}}$ is a (left) 
$\cal K$-principal
homogeneous space and the morphism ${\cal T}\lgr {\Bbb M}_p$ is smooth and $\cal K$-equivariant.
Thus, to obtain a log-smooth resolution of ${\cal S}_K(G,\mu)_{{\Bbb Z}_{(p)}}$ it will be
enough to construct a $\cal K$-equivariant log-smooth resolution of ${\Bbb M}_p$~: if
we have such a resolution $\widetilde {\Bbb M}_p$ , the morphism ${\cal K}\backslash ({\cal
T}\times_{{\Bbb M}_p}\widetilde {\Bbb M}_p)\lgr {\cal S}_K(G,\mu)_{{\Bbb Z}_{(p)}}$ will be a
log-smooth resolution. So, the theorem A derives from the following result:
\th T{HEOREM} B
\enonce
The local model ${\Bbb M}_p$ has a canonical $\cal K$-equivariant log-smooth
resolution.
\endth

Let us now describe the content of this article.

In section 2, we develope a variant of a reduction principle due to G. Faltings, which reduces
the problem of constructing resolutions of the local model to the one of resolving a scheme
$\mu$ defined by matrix equations.

In section 3, we use Lafforgue's work  on the compactifications of ${\rm
PGL}^{N+1}_r/{\rm PGL}_r$ ([L~1,2]) to solve this last problem.

In section 4, we give remarks concerning Faltings's original reduction principle and our
variant.
\vskip3mm
I want to thank A. J. de Jong, B. H. Gross, L. Illusie, L. Lafforgue, G. Laumon, M. Rapoport
and I. Vidal
for interesting discussions about this problem. I also want to thank G. Laumon for having
encouraged me to write this article.
\vskip7mm
{\bf 2\pointir A reduction}\vskip3mm
2.1\pointir In [F 2], G. Faltings introduces the scheme
$$\matrix{\mu^{r,N}=\{(A_0,\cdots,A_{N})\in {\goth gl}_r^{N+1}\mid &A_0A_1\cdots A_{N}&=&p.{\rm
Id}_r\hfill\cr
&A_1\cdots A_{N}A_0&=&p.{\rm Id}_r\hfill\cr
&&\vdots\cr
&A_{N}A_0\cdots A_{N-1}&=&p.{\rm Id}_r\}, } $$
endowed with a right action of ${\rm GL}_r^{N+1}$ defined by
$$(g_0,\cdots,g_{N}), (A_0,\cdots,A_{N})\mapsto
(g_{N}^{-1}A_0g_0\,,\,g_0^{-1}A_1g_1\,,\,\cdots,g_{N-1}A_{N}g_{N}).
$$
and proposes the following strategy to resolve the local models.

1)  Establish a reduction principle connecting ${\cal K}$-equivariant resolutions of ${\Bbb
M}_p$ to
${\rm GL}_r^{N+1}$-equivariant resolutions of $\mu^{r,N}$ (resp. with $r=g$ in the symplectic
case).

2) Construct such a ${\rm GL}_r^{N+1}$-equivariant resolutions of $\mu^{r,N}$.

In the present work, we will use a variant of Faltings's original reduction principle:
we shall replace ${\rm GL}_r$ by a parabolic subgroup $P$  of ${\rm GL}_n$ with Levi subgroup
${\rm GL}_r\times {\rm GL}_{n-r}$ (with $r=n-r=g$ in the symplectic case) and $\mu^{r,N}$ by an
analogous scheme $\mu^{P,N}$ also defined by matrix equations ({\it cf.} 2.2 for the unitary
case and 2.3 for the symplectic one). 
\vskip5mm
2.2\pointir The unitary case
\vskip2mm
Let $(d_i)_{1\leq i\leq {N+1}}$ be a $(N+1)$-uple of strictly positive integers with sum $n$.
Let $T$ be the matrix
$$\pmatrix{
0&{\rm Id}_{n-1}\cr
t&0
}\in{\goth gl}_n({\Bbb Z}[t]),$$
and ${\cal V}_i\  (0\leq i\leq N+1)$ be the lattice $T^{\,d_1+\cdots+d_i}{\Bbb
Z}[t]^n\subset{\Bbb Z}[t,t^{-1}]^n$. We denote by $\alpha_i\ 
(1\leq i\leq N+1)$ the natural inclusion
${\cal V}_i\inj {\cal V}_{i-1}$ and by $\alpha_0$ the composite inclusion ${\cal V}_0=\!={\cal
V}_{N+1}\inj{\cal V}_{N}$.
 \th D{EFINITION} 2.2.1. (cf. [R], [RZ]) 
\enonce The local model ${\Bbb M}_t\,({\rm
GL}_n\,,\,r\,,\,(d_i)_i)$ consists of
$(N+1)$-uples $(\omega_i\subset{\cal V}_i)_i$, where $$\omega _i \hbox{ is locally a direct
factor of rank
$r$ of ${\cal V}_i$}$$ and $$\alpha_i\,(\omega_i)\subset\omega_{i-1}\,,\forall i\in {\Bbb
Z}/(N+1){\Bbb Z}
$$

The local model ${\Bbb M}_p\,({\rm
GL}_n,r,(d_i)_i)$ is the restriction of ${\Bbb M}_t\,({\rm
GL}_n,r,(d_i)_i)$ along  the section $\overline p=\{t=p\}$ of the affine line ${\Bbb A}^1_t={\rm
Spec\,}{\Bbb Z}\,[t]$.
\endth
The functor ${\Bbb M}_t\,({\rm GL}_n,r,(d_i)_i)$ is obviously representable by a closed subscheme
in the product of grassmannian schemes $\prod_i{\rm Gr\,}(r,{\cal V}_i)$ ; the scheme ${\Bbb
M}_t\,({\rm GL}_n,r,(d_i)_i)[t^{-1}]$ is simply ${\rm Gr}(r,{\Bbb Z}[t,t^{-1}]^n)$.

Let ${\cal P}_t$ be the automorphism group of the system
$({\cal V}_i\,,\alpha_i)_i$ \ ({\it i.e.} ${\cal P}_t$ consists of $(N+1)$-uples
$(g_i)_i\in\prod_i{\rm GL}({\cal V}_i)$ such that $g_{i-1}\alpha_i=\alpha_i\, g_i\,,\ \forall
i\in {\Bbb Z}(N+1){\Bbb Z}$). The group scheme ${\cal P}_t$ acts obviously (on the left) 
on ${\Bbb M}_t\,({\rm GL}_n,r,(d_i)_i)$.
It is affine and smooth over ${\Bbb A}^1_t$ and ${{\cal P}_t}[t^{-1}]$ is naturally
identified with ${\rm GL}({\Bbb Z}[t,t^{-1}]^n)$. The group scheme $\cal K$ occuring in the
first section  is (in the unitary case) the restriction ${\cal P}_p$ of ${\cal P}_t$ along the
section
$\overline p$.

We will now define a certain rigidification of ${\Bbb M}_t\,({\rm GL}_n,r,(d_i)_i)$.
Let us denote by $V_i\ (i\in {\Bbb Z}(N+1){\Bbb Z})$ the vector bundle ${\Bbb A}^n$
and by 
$\Omega_i$ its direct factor ${\Bbb A}^r$.
 \th D{EFINITION} 2.2.2
\enonce The rigidified local model ${\Bbb M}_t^{\rm\, rig}\,({\rm
GL}_n,r,(d_i)_i)$ consists of $2(N+1)$-uples
$$((\omega_i\subset{\cal V}_i)_i\,,\ (\phi_i)_i),$$ where $$(\omega_i\subset{\cal
V}_i)_i \in {\Bbb M}_t\,({\rm GL}_n,r,(d_i)_i)$$ and $$\phi_i\ :\ (\omega_i\subset{\cal
V}_i)\lgr(\Omega_i\subset V_i)_{{\Bbb Z}[t]}$$ is an isomorphism of filtered vector bundles,
$\forall i\in {\Bbb Z}(N+1){\Bbb Z}$.
\endth
Let $P$ be the parabolic subgroup
$$\pmatrix{{\rm GL}_{r}&*\cr
0&{\rm GL}_{n-r}}$$
of ${\rm GL}_{n}$. The scheme ${\Bbb M}_t^{\rm\, rig}\,({\rm
GL}_n,r,(d_i)_i)$ is naturally a (right) $P^{N+1}$-principal homogeneous space over ${\Bbb M}
\,({\rm
GL}_n,r,(d_i)_i)$. It is also
endowed with  a (left) free action of
${\cal P}_t$ and both actions commute.
\vskip 2mm
2.2.3\pointir Let $\goth P$ be the Lie algebra of $P$. We will now re-write the algebraic space
${{\cal P}_t}\backslash{\Bbb M}_t^{\rm\, rig}\,({\rm
GL}_n,r,(d_i)_i)$
as a locally closed subscheme of $\goth P^{N+1}\times{\Bbb A}^1_t$.

Let $\goth P_{\rm rk\geq 1}$ be the open subscheme $\goth P-\{0\}$ of the affine space $\goth P$
and $\mu^{P,N}$ be the closed subscheme of $\goth P_{\rm rk\geq 1}^{N}\times{\Bbb A}^1_t$
defined by the matrix equations
$$\matrix{\Pi_0\Pi_1\hfill\cdots\hfill\Pi_{N}&=&t.{\rm Id}_n\cr
\Pi_{1}\hfill\cdots\hfill\Pi_{N}\Pi_0&=&t.{\rm Id}_n\cr
&\vdots&\cr
\Pi_{N}\Pi_0\cdots\Pi_{N-1}&=&t.{\rm Id}_n\,.\cr
}$$
The $(N+1)$-uple of matrices $$\Pi_j:=\phi^{-1}_{j-1}\alpha_j\phi_j\in\goth P_{\rm rk\geq
1}({\Bbb M}_t^{\rm\, rig}\,({\rm GL}_n,r,(d_i)_i))\qquad (j\in\{0,\cdots, N\})$$ clearly defines
a morphism
$${\Bbb M}_t^{\rm\, rig}\,({\rm
GL}_n,r,(d_i)_i)\lgr\mu^{P,N}$$
of ${\Bbb A}^1_t$-schemes, with image contained in the biggest open subscheme
$\mu^{P,N}_{(d_i)_i}$ of $\mu^{P,N}$, on which the minors of size $(n-d_i)$ of the matrix $\Pi_i$
are invertible, $\forall i\in {\Bbb Z}(N+1){\Bbb Z}$. Moreover, this morphism is
$P^{N+1}$-equivariant if we endow  $\mu^{P,N}$ with the action of $P^{N+1}$ defined by
$$(g_0,\cdots,g_N), (P_0,\cdots,P_N)\mapsto
(g_N^{-1}P_0g_0\,,\,g_0^{-1}P_1g_1\,,\,\cdots,g_{N-1}P_Ng_N).
$$
\th P{ROPOSITION} 2.2.4
\enonce The $\mu^{P,N}_{(d_i)_i}$-scheme ${\Bbb M}_t^{\rm\,
rig}\,({\rm GL}_n,r,(d_i)_i)$ is a ${\cal P}_t$-principal homogeneous space
\endth
{\it Proof} ({\it cf.} [RZ], appendix A): It is clearly sufficient  to show that locally for the
Zariski topology on
$\mu^{P,N}_{(d_i)_i}$ the morphism $${\Bbb M}_t^{\rm\, rig}\,({\rm
GL}_n,r,(d_i)_i)\lgr\mu_{(d_i)_i}^{P,N}$$ has a section. So, let $R$ be a local ring of
$\mu^{P,N}_{(d_i)_i}$, with maximal ideal $\cal M$ and residue field $k$. The case $t\in
R^\times $ is trivial, and we assume $t\in {\cal M}$. The matrix $\Pi_i\otimes k$ is of rank
$\geq (n-d_i)$, hence the 
$k$ vectorspace
$(V_{i-1}\otimes k)/\Pi_i(V_{i}\otimes k)$ is of rank $d'_i\leq d_i$. Let $(e_j[i])_{1\leq
j\leq d'_i}$ be a basis of this $k$-vectorspace and let $(\widetilde e_j[i])_{1\leq
j\leq d'_i}$ be a lifting of this family to a family of elements in $V_i\otimes R$. The family
$$(\widetilde e_j[i])_{1\leq
j\leq d'_i}\cup(\Pi_{i+1}\widetilde e_j[i+1])_{1\leq
j\leq d'_{i+1}}\cup\cdots\cup(\Pi_{i+1}\cdots\Pi_{i+N}\widetilde e_j[i+N])_{1\leq
j\leq d'_{i+N}} $$
(where $i\,$ is understood as an element of ${\Bbb Z}(N+1){\Bbb Z}$) generates $V_i\otimes k$,
hence we have $d'_j=d_j\,,\  \forall j$. Let $U_i$ be the $R$-submodule of $V_i\otimes R$
spanned by the family $(\widetilde e_j[i])_{1\leq
j\leq d'_i}$. This family obviously defines a morphism $R^{d_i}\lgr U_i$. Let us
consider the morphism
$$\psi_i\ :\ R^{d_i}\oplus R^{d_{i+1}}\oplus\cdots\oplus R^{d_{i+N}}\lgr V_i\otimes R
$$ induced by $R^{d_j}\lgr U_j \ (i\leq j\leq i+N)$ and $\Pi_{i+1}\cdots\Pi_j\ :\
U_j\lgr V_i\otimes R\  (i+1\leq j \leq i+N)$. It is surjective by Nakayama's  Lemma, and so
it must be an isomorphism. The composite $$U_i=\!=R^{d_{i+N}}\lgr
V_{i-1}=\!=V_{i+N}$$ of its last factor with
$\Pi_i\ :\ V_i\otimes R\lgr V_{i-1}\otimes R$ is the composite
$(\Pi_i\cdots\Pi_{i+N})_{\vert U_i}=(t.{\rm Id}_n)_{\vert U_i}$. Hence, the matrix (in the
canonical bases) of $\psi_{i-1}^{-1}\Pi_i\psi_i$ is simply $T^{d_i}$ and the $2(N+1)$-uple
$((\Pi_i)_i,(\psi_i^{-1})_i)$ defines a section of ${\Bbb M}_t^{\rm\, rig}\,({\rm
GL}_n,r,(d_i)_i)_R\lgr(\mu_{(d_i)_i}^{P,N})_R$ .\quad $\square$
\vskip2mm
2.2.5\pointir We will now state our variant of  G. Faltings's reduction principle.

Let $X$ be a ${\Bbb A}^1_t$-scheme. We assume that $X[t^{-1}]$ is smooth over ${\Bbb
A}^1_t[t^{-1}]={\rm Spec\, {\Bbb Z}\,}[t,t^{-1}]$.
\th D{EFINITION} 2.2.6
\enonce A toro\"\i dal resolution of $X$ is a proper morphism
of ${\Bbb A}^1_t$-schemes $r\ :\  \widetilde X\lgr X$ such that $r[t^{-1}]$ is an
isomorphism and such that there exists a torus $T$, a character $\chi$ of $T$, a toro\"\i dal
embedding
$T\inj
\overline T$ (cf. {\rm [KKMS]})  and  a diagram
$$\diagram{
&\ \widetilde X^\vee
\cr\noalign{\vskip 5mm}
\widetilde X
&&\overline T
\cr
\fleche(7,9)\dir(-1,-1)\long{4}
\fleche(11,9)\dir(1,-1)\long{4}
}$$
of ${\Bbb A}^1_t$-schemes satisfying the following conditions

\item{--} the restriction to $T$ of the morphism $\overline T\lgr{\Bbb A}^1_t$ is the composite
$T\buildrel \chi\over\lgr{\Bbb G}_m\inj{\Bbb A}^1_t$

\item{--} the left morphism is a principal homogeneous space under a torus $T'$ endowed with a
morphism to the kernel
$T_\chi$ of
$\chi$

\item{--} the right morphism is smooth and $T'$-equivariant.

We say this toro\"\i dal resolution is strongly log-smooth if moreover the group $T_\chi$ is
a torus.

\endth
 \th R{EDUCTION  PRINCIPLE} 2.2.7 
\enonce Assume that we have 
a $P^{N+1}$-equivariant strongly log-smooth toro\"\i dal resolution 
$\widetilde\mu^{P,N}$ of $\mu^{P,N}$ whose source can be put in a
$P^{N+1}\times T'$-equivariant diagram {\rm (2.2.6)}, where  the action of $P^{N+1}\times T'$ on
$\overline T$ is by the second factor. The morphism
$$\widetilde{\Bbb M}_t\,({\rm
GL}_n,r,(d_i)_i):=({\Bbb M}_t^{\rm\, rig}\,({\rm
GL}_n,r,(d_i)_i)\times_{\mu^{P,N}}\widetilde\mu^{P,N})/P^{N+1}\lgr{\Bbb M}_t\,({\rm
GL}_n,r,(d_i)_i) $$
is then a ${\cal P}_t$-equivariant strongly log-smooth toro\"\i dal resolution of ${\Bbb
M}_t\,({\rm GL}_n,r,(d_i)_i)$.
\endth
{\it Proof}: This morphism is clearly proper, and an isomorphism over ${\Bbb A}^1_t[t^{-1}]$.
As a diagram (2.2.6) for $\widetilde{\Bbb M}_t\,({\rm
GL}_n,r,(d_i)_i)$ we can take $$(({\Bbb M}_t^{\rm\, rig}\,({\rm
GL}_n,r,(d_i)_i)\times_{\mu^{P,N}}[\hbox{diagram (2.2.7) for $\widetilde \mu^{P,N}$}])/P^{N+1}
.$$ The only non-formal point to check is that the morphism $({\Bbb M}_t^{\rm\, rig}\,({\rm
GL}_n,r,(d_i)_i)\times_{\mu^{P,N}}\widetilde\mu^{P,N})/P^{N+1}\lgr \overline T_t$ is smooth.
This follows from  proposition (2.2.4). $\square$
\vskip2mm
Taking for granted the existence of a $P^{N+1}$-equivariant strongly log-smooth toro\"\i dal
resolution 
$\widetilde\mu^{P,N}$ of $\mu^{P,N}$  ({\it cf}. section
3), the (unitary case of the) theorem announced in the first section will be a consequence of
this reduction principle and of the following proposition (which also explains the terminology
2.2.6)
  \th P{ROPOSITION} 2.2.8
\enonce

 a) Let $T\inj \overline T$ be a torus embedding, and $\overline \chi\ :\ \overline T\lgr
{\Bbb A}^1_t$ with restriction to $T$ given by a character $\chi$ of $T$. If we assume
that the kernel
$T_\chi$ of $\chi$ is a torus and that $\overline T$ and
${\Bbb A}^1_t$ are respectively endowed with the canonical log-structures associated with the
torus embeddings
$T\inj \overline T$ and
${\Bbb G}_m\inj{\Bbb A}^1_t$, the morphism
$\overline T_t$ is log-smooth.

b) As a particular case, under the conditions of (a),
the restriction $\overline T_p$ of $\overline T_t$ over the section $\overline p=\{t=p\}$ of
${\Bbb A}^1_t$ is log-smooth

c) Let
$$\diagram{
& Y
\cr\noalign{\vskip 6mm}
 X
&\lgr& Z
\cr
\fleche(7,8)\dir(-1,-1)\long{4}
\fleche(11,8)\dir(1,-1)\long{4}
}$$
be a commutative triangle of morphisms of  fine log-schemes. If we assume that the left morphism
is smooth and surjective, that the log-structure on $Y$ is the inverse image of the
log-structure on
$X$, and that $Y$ is log-smooth over $Z$, then $X$ is log-smooth over $Z$
\endth
{\it Proof}: a) Let $X^*(T)$ be the character group of the torus $T$. The kernel of $\chi$ is a
torus, and so the cokernel of the morphism ${\Bbb Z}\chi\lgr X^*(T)$ is torsion-free. Hence
the morphism $\overline T_t$ is log-smooth.

b) It is a direct consequence of (a).

c) (after I. Vidal) Log-smoothness is a local property for the (classical) \'etale topology,
hence we can asssume that $Y$ is an affine space ${\Bbb A}^d_X$ over $X$. In this case, we will
use Kato's infinitesimal lifting criterion. Let
$S_0\inj S$ be a strict nilimmersion ; every morphism of log-schemes $S_0\lgr X$ can be lifted
to a morphism of log-schemes $S_0\lgr {\Bbb A}^d_X$ and the log-smoothness of ${\Bbb
A}^d_X=Y\lgr Z$ implies that  $S_0\lgr {\Bbb A}^d_X$ can be lifted to $S\lgr {\Bbb A}^d_X$. The 
projection $S\lgr X$ of this last morphism is obviously a  lifting of $S_0\lgr X$, and so
the proposition is proved. $\square$
\vskip 5mm
2.3\pointir The symplectic case
\vskip2mm
Let $(d_i)_{1\leq i\leq N}$ be a $N$-uple of strictly positive integers with sum $g$ and
let $d_{N+1}=g$. The $(N+1)$-uple $(d_i)_{1\leq i\leq N+1}$ is a partition of $2g$ and the local
model in the symplectic case, ${\Bbb M}_t\,({\rm
Sp}_{2g},(d_i)_i)$, is the following closed subscheme of ${\Bbb M}_t\,({\rm
GL}_{2g}\,,\,g,(d_i)_{1\leq i\leq N+1})$.

Let $\Delta$ be the matrix associated with the permutation $(g,g-1,\cdots,1)$ and $J\in {\rm
GL}_{2g}({\Bbb Z})$ be the matrix
$$\pmatrix{0&-\Delta\cr
\Delta&\hphantom{-}0}$$
The matrix $J$ defines on ${\Bbb Z}[t,t^{-1}]$ a nondegenerate symplectic form $<\ ,\ >$. Its
restriction $<\ ,\ >_0$ to ${\cal V}_0$ is nondegenerate and the restriction $<\ ,\ >_g$ of
$t^{-1}<\ ,\ >$ to ${\cal V}_g$ has values in ${\Bbb Z}[t]$ and is also nondegenerate.
\th D{EFINITION} 2.3.1. (cf. [dJ 2], [RZ])
\enonce
 The local model ${\Bbb M}_t\,({\rm
Sp}_{2g},(d_i)_i)$ consists of
$(N+1)$-uples $$(\omega_i\subset{\cal V}_i)_i\in {\Bbb M}_t\,({\rm
GL}_{2g}\,,\,g,(d_i)_{1\leq i\leq N+1})$$ where $\omega_0$ and  $\omega_N$ are respectively
totally isotropic with respect to $<\ ,\ >_0$  and $<\ ,\ >_N$).

The local model ${\Bbb M}_p\,({\rm
Sp}_{2g},(d_i)_i)$ is the restriction of ${\Bbb M}_t\,({\rm
Sp}_{2g},(d_i)_i)$ along  the section $\overline p=\{t=p\}$ of the affine line 
${\Bbb A}^1_t={\rm
Spec\,}{\Bbb Z}\,[t]$.
\endth

The scheme
${\Bbb M}_t\,({\rm Sp}_{2g},(d_i)_i)[t^{-1}]$ is simply the lagrangian grassmannian of (${\Bbb
Z}[t,t^{-1}]^{2g}\,,\break <\ ,\ >$).

Let ${\cal Q}_t$ be the automorphism group of the system
($({\cal V}_i\,,\alpha_i)_i,<\ ,\ >_0\,,\,<\ ,\ >_N$) \ ({\it i.e.} ${\cal Q}_t$ is the closed
subgroup of ${\cal P}_t$ defined by the conditions $\gamma_0\in {\rm Sp}\,({\cal V}_0\,,\,<\ ,\
>_0$ and $\gamma_N\in {\rm Sp}\,({\cal V}_N\,,\,<\ ,\
>_N$)). The group scheme
${\cal Q}_t$ acts obviously (on the left)  on
${\Bbb M}_t\,({\rm Sp}_{2g},(d_i)_i)$. It is affine and smooth over ${\Bbb A}^1_t$ and ${{\cal
Q}_t}[t^{-1}]$ is naturally identified with ${\rm Sp}\,({\Bbb Z}[t,t^{-1}]^{2g},<\ ,\ >)$. The
group scheme $\cal J$ occuring in the first section  is (in the symplectic case) the restriction
${\cal Q}_p$ of
${\cal Q}_t$ along the section
$\overline p$.

We will now define a certain rigidification of ${\Bbb M}_t\,({\rm Sp}_{2g}\,,\,(d_i)_i)$.
 \th D{EFINITION} 2.3.2
\enonce
 The rigidified local model ${\Bbb M}_t^{\rm\, rig}\,({\rm
Sp}_{2g}\,,\,(d_i)_i)$ is the fibre product 
$${\Bbb M}_t\,({\rm Sp}_{2g}\,,\,(d_i)_i)\times_{{\Bbb M}_t\,({\rm
GL}_{2g}\,,\,g\,,\,(d_i)_i)}{\Bbb M}_t^{\rm rig}\,({\rm GL}_{2g}\,,\,g\,,\,(d_i)_i). $$
\endth
The scheme ${\Bbb M}_t^{\rm\, rig}\,({\rm
Sp}_{2g}\,,\,(d_i)_i)$ is naturally a (right) $P^{N+1}$-principal homogeneous space over 
${\Bbb M}_t\,({\rm
Sp}_{2g}\,,\,(d_i)_i)$. It is also
endowed with  a (left) free action of
${\cal Q}_t$ and both actions commute.
\vskip 2mm
2.3.3\pointir We will now re-write the algebraic space
${{\cal Q}_t}\backslash{\Bbb M}_t^{\rm\, rig}\,({\rm
Sp}_{2g}\,,\,(d_i)_i)$.

Let $\Sigma^{g,N}$ be the scheme consisting of triples $((\Pi_i)_{0\leq i\leq N}\,,\, (\ ,\
)_0\,,\,(\ ,\ )_N)$, where 

\item{--} $(\Pi_i)_{0\leq i\leq N}\in \mu^{P,N}$

\item{--} $(\ ,\
)_0$ (resp. $(\ ,\ )_N)$) is a nondegenerate symplectic pairing on $V_0$ (resp. $V_N$) such
that $\Omega_0$ (resp. $\Omega_N$) is totally isotropic

\item{--} $(\Pi_0\cdots\Pi_{N-1}\, x,y)_N=(x,\Pi_{N}y)_0\,,\ \forall x\in V_0\,,\ y \in V_N$

\noindent and let $\Sigma^{g,N}_{(d_i)_i}$ be the fibre product
$\Sigma^{g,N}\times_{\mu^{P,N}}\mu^{P,N}_{(d_i)_i}$. 
Using (2.2.4) and a variant of ([dJ~2], proposition 3.6)
(or of  [RZ], theorem 3.16)  we obtain the following statement. 
  \th P{ROPOSITION} 2.3.4
\enonce
 The $\Sigma^{P,N}_{(d_i)_i}$-scheme ${\Bbb M}_t^{\rm\,
rig}\,({\rm Sp}_{2g},(d_i)_i)$ is a ${\cal Q}_t$-principal homogeneous space.
\endth
{\it Proof}: the proof of [ibid.] carries over easily to our context (``$p$ replaced by $t$").
Using (2.2.4), the system  $((V_i)_i, (\Pi_i)_i, (\ ,\ )_0, (\ ,\ )_N)$ is a system of type II
(in the terminology of [dJ 2]) and so the systems 
$((V_i)_i, (\Pi_i)_i, (\ ,\ )_0, (\ ,\ )_N)$ and $(({\cal V}_i)_i\,,\, (\alpha_i)_i, <\ ,\ >_0,
<\ ,\ >_N)$ are (Zariski-)locally isomorphic  on $\Sigma^{P,N}_{(d_i)_i}$. $\square$
\vskip2mm
2.3.5\pointir The morphism $\Sigma^{g,N}_{(d_i)_i}\lgr\mu^{r,N}$ is not smooth. However, as
Faltings remarked in his analogous situation, one can turn out this difficulty by using toro\"\i
dal resolutions with the following additional property. 

Let us write
$$\Pi_{N}=\pmatrix{A_{N}&B_{N}\cr0&C_{N}}$$
and 
$$\Pi_0\cdots\Pi_{N-1}=\pmatrix{A'_{N}&B'_{N}\cr0&C'_{N}}$$ 
and let $\mu^{P,N}_{\rm bl.}$ be the scheme obtained by blowing up $\mu^{P,N}$ along the ideal
sheaves spanned by the minors of size $i$ of $A_{N}$, $\forall\, 1\leq i\leq g-1$ and also along
the ideal sheaves spanned by the minors of size $i$ of $A'_{N}$, $\forall\, 1\leq i\leq g-1$.

We shall only consider resolutions of $\mu^{P,N}$ factoring through the morphism 
$\mu^{P,N}_{\rm bl.}\lgr\mu^{P,N}$.
The following proposition, essentially due to Faltings ([F 2]), is the crucial step of the
reduction principle for these resolutions.

For any ${\Bbb Z}[t]$-scheme $X$ let us denote by $X^+$ the closed subsheme of
$X$ defined by the ideal sheaf $\hbox{\it Tors}_t(X)$ consisting of sections of ${\cal O}_{X}$
killed by a power of
$t$.
\th P{ROPOSITION} 2.3.6
\enonce
 The morphism
$$(\Sigma^{g,N}_{(d_i)_i}\times_{\mu^{P,N}}\mu^{P,N}_{\rm bl.})^+\lgr (\mu^{P,N}_{\rm bl.})^+$$
is smooth.
\endth
{\it Proof}: Let $(\Sigma^{g,N})^0_B$ (resp.$(\Sigma^{g,N})^0_{B'}$) be the  biggest  open
subscheme of
$\Sigma^{g,N}_{(d_i)_i}$ where the matrix $B_{N}$ (resp.  $B'_{N}$) is  invertible. It is
easily seen that the translates of $(\Sigma^{g,N})^0_B$ (resp.$(\Sigma^{g,N})^0_{B'}$) under
$P^{N+1}$ cover $\Sigma^{g,N}_{(d_i)_i}$. Hence their intersection $(\Sigma^{g,N})^0$ is also such
that its translates under $P^{N+1}$ cover $\Sigma^{g,N}_{(d_i)_i}$. The morphism
$(\Sigma^{g,N}_{(d_i)_i}\times_{\mu^{P,N}}\mu^{P,N}_{\rm bl.})^+\lgr (\mu^{P,N}_{\rm bl.})^+$
is $P^{N+1}$-equivariant and so it suffices to prove that the morphism 
$((\Sigma^{g,N})^0\times_{\mu^{P,N}}\mu^{P,N}_{\rm bl.})^+\lgr (\mu^{P,N}_{\rm bl.})^+$ is
smooth.

Changing the bases of $V_0$ and $V_N$ we can assume that $B_{N}=B'_{N}={\rm Id}_g$ and
(because the ideal sheaves spanned by the minors of size $i$ of $A_{N}$ (resp. $A'_{N}$)
have been blown up, $\forall \,1\leq  i\leq g-1$) that the matrix $A_{N}$ (resp. $A'_{N})$ is
the diagonal matrix 
$${\rm diag\,}(a_0,a_0a_1,\cdots,a_0\cdots a_{g-1})$$ 
(resp. $${\rm diag\,}(a_1\cdots a_g,a_2\cdots a_g,\cdots,a_g),$$ 
with $a_0a_1\cdots a_g=t$). In the new bases we have
$$\Pi_{N}=\pmatrix{A_{N}&{\rm Id}_g\cr
0&-A'_{N}}\quad{\rm and}\quad \Pi'_{N}=\pmatrix{A'_{N}&{\rm Id}_g\cr
0&-A_{N}}$$
(recall that $\Pi_{N}\Pi'_{N}=\Pi'_{N}\Pi_{N}=t.{\rm Id}_{2g}$). Let 
$$J[0]=\pmatrix{0&J_1[0]\cr
J_2[0]&J_3[0]}\quad{\rm and}\quad J[N]=\pmatrix{0&J_1[N]\cr
J_2[N]&J_3[N]}$$
be the matrices defining the symplectic pairings on $V_0$ and $V_N$ in the new bases (the upper
left terms vanish because $\Omega_0$ and $\Omega_N$ are totally isotropic). These matrices
satisfy the following relations
$$\displaylines{{}^tJ_1[0]=-J_2[0]\,,\ {}^tJ_3[0]=-J_3[0]\,,\  {}^tJ_1[N]=-J_2[N]\,,\
{}^tJ_3[N]=-J_3[N]\,,\cr
A'_{N}J_1[N]+J_1[0]A'_{N}=0\,,\cr
 J_1[N]A_{N}+A_{N}J_1[0]=0\,,\cr
J_1[N]-A_{N}J_3[N]+{}^tJ_1[0]+J_3[0]A'_{N}=0 } $$
Using the first line, we eliminate the coefficients of $J_2[0]$ and of $J_2[N]$. Killing the
$t$-torsion, these relations become
$$\displaylines{J_1[0]^i_j+a_i\cdots a_{j-1} J_1[N]^i_j=0\,,\ \forall \,1\leq i\leq j\leq g\cr
J_1[N]^i_j+a_j\cdots a_{i-1}J_1[0]^i_j=0\,,\ \forall \,1\leq j\leq i\leq g\cr
J_1[N]^i_j+J_1[0]^j_i-a_0\cdots a_{i-1}J_3[N]^i_j+a_j\cdots a_g J_3[0]^i_j=0\,,\ \forall
\,1\leq i\leq j\leq g }$$
Using the two first lines, we eliminate the coefficients $J_1[0]\ (i\leq j)$ and $J_1[N]\
(i>j)$. Using the last one, we eliminate $J_1[0]\ (i>j)$. The other coordinates $J_1[N]\ (i\leq
j)\,,\ J_3[0]\ (i<j)\,,\ J_3[N]\ (i<j)$ can be fixed freely, and so we see that 
$((\Sigma^{g,N})^0\times_{\mu^{P,N}}\mu^{P,N}_{\rm bl.})^+$ is an affine space of dimension
$g(3g-1)/2$ over $(\mu^{P,N}_{\rm bl.})^+$. This achieves the proof of  (2.3.6). $\square$

Using the proposition (2.3.6) we obviously obtain the following variant of Faltings's reduction
principle.
 \th R{EDUCTION  PRINCIPLE} 2.3.7 
\enonce
 Assume that we have 
a $P^{N+1}$-equivariant toro\"\i dal resolution 
$\widetilde\mu^{P,N}\lgr\mu^{P,N}_{\rm bl.}\lgr\mu^{P,N}$ whose source can be put in a
$P^{N+1}\times T'$-equivariant diagram {\rm (2.2.6)}, where  the action of $P^{N+1}\times T'$ on
$\overline T_t$ is by the second factor. The morphism
$$\widetilde{\Bbb M}_t\,({\rm
Sp}_{2g},(d_i)_i):=({\Bbb M}_t^{\rm\, rig}\,({\rm
Sp}_{2g},(d_i)_i)\times_{\mu^{P,N}}\widetilde\mu^{P,N})/P^{N+1}\lgr{\Bbb M}_t\,({\rm
Sp}_{2g},(d_i)_i) $$
is then a ${\cal Q}_t$-equivariant toro\"\i dal resolution of ${\Bbb
M}_t\,({\rm Sp}_{2g},(d_i)_i)$.$\square$ 
\endth
Taking for granted the existence of a $P^{N+1}$-equivariant  strongly log-smooth toro\"\i dal
resolution 
$\widetilde\mu^{P,N}$ of $\mu^{P,N}$ with this property ({\it cf}. section
3), the (symplectic case of the) theorem announced in the first section will be a consequence of
this reduction principle and of the  proposition (2.2.8).
\vskip 7mm
{\bf 3\pointir ``Compactification of the multiplication law" for parabolic subgroups of the
general linear group }
\vskip 3mm
We will now construct a $P^{N+1}$ -equivariant   strongly log-smooth toro\"\i dal resolution
of $\mu^{P,N}$.  Our construction follows quite closely Lafforgue's  compactification ([L 1,2])
of
${\rm PGL}_r^{N+1}/{\rm PGL}_r$.
\vskip 2mm 
3.1\pointir Let us first briefly recall Lafforgue's construction.
 
Let $P_{\!\! ad}$ be the quotient of $P$ by its center ${\Bbb G}_m$.
Lafforgue constructs a
$P_{\!\! ad}^{N+1}$ -equivariant open immersion 
$$P_{\!\! ad}^{N+1}/ P_{\!\! ad}\hookrightarrow \overline \Omega^{P,N}$$
as a quotient of an open immersion
$$ P^{N+1}\times {\Bbb G}_m^{S^{P,N}}/( P\times
{\Bbb G}_m^{N+1})\hookrightarrow \Omega^{P,N}$$ 
by a torus ${\Bbb G}_m^{S^{P,N}}/{\Bbb G}_m$ acting freely on both sides (see [L 1] for the
special case $P={\rm GL}_n$ and [L 2] for  arbitrary
parabolic subgroups of ${\rm GL}_n$). Moreover,
$\Omega^{P,N}$ is endowed with a
${\Bbb G}_m^{S^{P,N}}/{\Bbb G}_m$-equivariant morphism
$\Omega^{P,N}\longrightarrow {\cal A}^{P,N}$ to a toric normal variety ${\cal A}^{P,N}$
with torus ${\cal A}_{\emptyset}^{P,N}={\Bbb G}_m^{S^{P,N}}/{\Bbb G}_m^{N+1}$ (we shall
recall in the next paragraph the definition of the tori and of the maps between them).

The scheme $\overline\Omega^{P,N}$ is a ``toro\"\i dal equivariant compactification" of
$P_{\!\! ad}^{N+1}/ P_{\!\! ad}$ in the following sense ({\it cf.} [L 2],Th\'eor\`emes 5 et 6):
 the scheme $\overline\Omega^{P,N}$ is projective over ${\rm Spec}\,{\Bbb Z}$
and
 the morphism $\Omega^{P,N}\longrightarrow {\cal A}^{P,N}$ is smooth.

The scheme $\overline\Omega^{P,N}$ is a compactification of the multiplication law for $P_{\!\!
ad}$ in the following sense ([L~1, L 2]). For every map
$\{0,\cdots,M\}\longrightarrow\{0,\cdots ,N\}$, Lafforgue constructs commutative diagrams

$$\def\normalbaselines{\baselineskip20pt
  \lineskip3pt \lineskiplimit3pt }

\def\mapdown#1{\Big\downarrow
  \rlap{$\vcenter{\hbox{$\scriptstyle#1$}}$}}
\matrix{
     P^{N+1}\times {\Bbb G}_m^{S^{P,N}}/( P\times
{\Bbb G}_m^{N+1})&\hookrightarrow&\Omega^{P,N}&\cr
  \mapdown{}&&\,\mapdown{}\hfill\cr
     P^{M+1}\times {\Bbb G}_m^{S^{P,N}}/( P\times
{\Bbb G}_m^{M+1})&\hookrightarrow&
    \Omega^{P,M}&\cr
  }\qquad{\rm and}\qquad
\matrix{
     P_{\!\! ad}^{N+1}/ P_{\!\! ad}&\hookrightarrow&\overline\Omega^{P,N}&\cr
  \mapdown{}&&\,\mapdown{}\hfill\cr
     P_{\!\! ad}^{M+1}/ P_{\!\! ad}&\hookrightarrow&
    \overline\Omega^{P,M}&\cr
  }$$
\lineskip 3mm
(the diagram in the right hand side is of course a quotient of the one in the left hand side).

If as a particular case we consider the family of injections $$\{i,i+1\}\hookrightarrow
\{0,\cdots,N\}\quad (i\in {\Bbb Z}/(N+1){\Bbb Z}=\{0,\cdots,N\})$$ we obtain a morphism
$\overline\Omega^{P,N}\longrightarrow(\overline\Omega^{P,2})^{N+1}$ which is a
compactification of the closed immersion
$$\matrix{ P^{N+1}/ P\quad\ \hookrightarrow &(P^2/
P )^{N+1} &=\!=\  P^{N+1}\hfill\cr
(h_0,\cdots,h_{N})\hfill&\longmapsto &\qquad(h_0h_1^{-1},\cdots,h_{N_1}h_N^{-1},h_Nh_0^{-1})
}$$
with image $\{(g_0,\cdots, g_N)/g_0\cdots g_N=1\}$. So, the morphism
$\overline\Omega^{P,N}\longrightarrow(\overline\Omega^{P,2})^{N+1}$ is a compactification
of the morphism giving the product of $N$ elements in $ P$.
\vskip 2mm 
3.2\pointir In this paragraph, we shall modify the construction of $\overline\Omega^{P,N}$ to
keep track of the center of $P$. More precisely, we shall 
define a subtorus ${\cal T}^{P,N}$ of ${\Bbb
G}_m^{S^{P,N}}/\,{\Bbb G}_m$ such that the quotient $${\cal T}^{P,N}\backslash
P^{N+1}\times {\Bbb G}_m^{S^{P,N}}/( P\times {\Bbb G}_m^{N+1})$$ is 
${\Bbb G}_m\times P^{N+1}/ P\,$ and so the quotient ${\cal
T}^{P,N}\backslash\Omega^{P,N}$ will be a ``partial compactification" of ${\Bbb
G}_m\times P^{N+1}/ P$. We will also define a morphism 
$${\cal T}^{P,N}\backslash\Omega^{P,N}\longrightarrow {\goth
P}_{{\rm rk}\geq 1}^{N+1}\times{\Bbb A}^1_t
$$ factoring through the closed immersion $\mu^{P,N}\hookrightarrow {\goth
P}_{{\rm rk}\geq 1}^{N+1}\times{\Bbb A}^1_t$. 

The (projective) morphism
$$\widetilde\mu^{P,N}={\cal T}^{P,N}\backslash\Omega^{P,N}\longrightarrow {\mu}^{P,N}$$ thus
obtained will be the desired resolution.

\vskip 2mm
Let us recall the following definitions and facts ([L 2], 1.b).

$\bullet$ The finite set $S^{P,N}$ appearing in the beginning of the previous paragraph is
$$\{(i_{\alpha,1},i_{\alpha,2})_{0\leq\alpha\leq N}\in {\Bbb N}^{2(N+1)}\mid\sum_\alpha
\vert i\vert_\alpha=n\  {\rm and}\  \sum_\alpha
i_{\alpha,1}\geq r\}$$
(where, for simplicity, $\vert i\vert_\alpha$ denotes $i_{\alpha,1}+i_{\alpha,2}$).
It is endowed with the following partial order
$$i\leq j\buildrel\hbox{\fiverm def}\over\iff \vert i\vert_\alpha=\vert j\vert_\alpha \ {\rm
and}\ i_{\alpha, 1}\leq j_{\alpha, 1}\,,\ \forall \ 0\leq\alpha\leq N\ .$$

$\bullet$ The map $ P\times {\Bbb G}_m^{N+1}\hookrightarrow 
P^{N+1}\times {\Bbb G}_m^{S^{P,N}}$ is
$$(g,\lambda_0,\cdots,\lambda_N)\mapsto (g,\lambda_1 g,\cdots,\lambda_N g,(\lambda_0
\,{\rm det\,}g^{-1}\lambda_1^{-\vert i\vert_1}\cdots\lambda_N^{-\vert i\vert_N})_{i \in
S^{P,N}})\,.$$

$\bullet$ The map $
P^{N+1}\times {\Bbb G}_m^{S^{P,N}}/( P\times {\Bbb
G}_m^{N+1})\longrightarrow
P^{2}\times {\Bbb G}_m^{S^{P,2}}/( P\times {\Bbb G}_m^{2})$ associated
with an injection $\{j,j+1\}\hookrightarrow \{0,\cdots,N\}$  is simply
$$(g_0,\cdots,g_N,(\lambda_j)_j)\mapsto
(g_i,g_{i+1},(\lambda_{(0,\cdots,0,j_i,j_{i+1},0\cdots,0})_{\vert j\vert_i+\vert j\vert_{i+1}=n})
$$

$\bullet$ Let $E={\Bbb A}^{n(N+1)}$ be the trivial vector bundle of rank $n(N+1)$, equipped
with the obvious gradation $$E=\bigoplus_{0\leq\alpha\leq N}E_{\alpha}\qquad (E_{\alpha}\simeq
{\Bbb A}^n),$$ compatible filtration
$$({\Bbb A}^r=F_{\alpha}=E_{\alpha,1}\subset E_{\alpha,2}=E_{\alpha})$$ and action of $
P^{N+1}$. Let us recall the following notations from [L 2]

\item{--} One has $\Lambda^nE=\bigoplus_{\vert i\vert} \Lambda^{\vert i\vert}E$. The sum is over
$(N+1)$-uples
$$\vert i\vert=(\vert i\vert_0,\cdots,\vert i\vert_N)\in{\Bbb N}^{N+1}\ {\rm with}\ \vert
i\vert_1+\cdots \vert i\vert_N=n\,,$$ and $\Lambda^{\vert i\vert}E$ is
$\Lambda^{{\vert i\vert}_0}E_0\otimes\cdots\otimes\Lambda^{\vert i\vert_N}E_N$.

\item{--} For $i\in S^{P,N}$, $\Lambda^iE$ denotes the sub-vector bundle
$\bigotimes_\alpha(\Lambda^{i_{\alpha,1}}F_\alpha\wedge \Lambda^{i_{\alpha, 2}}E_\alpha)$ of
$\Lambda^{\vert i
\vert}E$.

$\bullet$ The scheme $\Omega^{P,N}$ inherits from its construction as a Zariski closure in the
product 
$${\Bbb G}_m\backslash \prod_{i\in S^{P,N}}[(\Lambda^{\vert
i\vert}E/\sum_{j>i}\Lambda^jE)-\{0\}]\times{\cal A}^{P,N}$$
a morphism
$$\Omega^{P,N}\longrightarrow {\Bbb G}_m\backslash \prod_{i\in S^{P,N}}[(\Lambda^{\vert
i\vert}E/\sum_{j>i}\Lambda^jE)-\{0\}]$$
This morphism is generically an immersion. Moreover it is projective 
([L 1,2], th\'eor\`eme~5).

$\bullet$ The composite of the open immersion $ P^{N+1}\times {\Bbb
G}_m^{S^{P,N}}/( P\times {\Bbb G}_m^{N+1})\hookrightarrow \Omega^{P,N}$ with
this morphism is
$$(g_0,\cdots,g_N,(\lambda_i)_i)\mapsto(\lambda_i\,\Lambda^{{\vert
i\vert}_0}({}^tg_0)\wedge\cdots\wedge\Lambda^{{\vert
i\vert}_N}({}^tg_N))_i$$
where $\Lambda^{{\vert
i\vert}_0}({}^tg_0)\,\wedge\,\cdots\,\wedge\,\Lambda^{{\vert
i\vert}_N}({}^tg_N)$ is understood as a morphism 
$$\Lambda^{{\vert i\vert}_0}E_0^\vee\,\otimes\,\cdots\,\otimes\,\Lambda^{\vert
i\vert_N}E_N^\vee\,\longrightarrow\Lambda^{{\vert
i\vert}_0}E_0^\vee\,\wedge\,\cdots\,\wedge\,\Lambda^{\vert i\vert_N}E_N^\vee\,={\Bbb
A}^1$$ (the non-invariance of this last identification disappears in the quotient), identified to
an element of
$\Lambda^{\vert i\vert}E$ by biduality.

$\bullet$ The composite  of the open immersion $ P^{N+1}\times {\Bbb
G}_m^{S^{P,N}}/( P\times {\Bbb G}_m^{N+1})\hookrightarrow \Omega^{P,N}$ with
the morphism $\Omega^{P,N}\longrightarrow{\cal A}^{P,N}$ is the obvious morphism $
P^{N+1}\times {\Bbb G}_m^{S^{P,N}}/( P\times {\Bbb
G}_m^{N+1})\longrightarrow {\Bbb G}_m^{S^{P,N}}/{\Bbb G}_m^{N+1}=\!={\cal
A}_{\emptyset}^{P,N}\hookrightarrow {\cal A}^{P,N}$. 

$\bullet$ One has a morphism $\Omega^{P,N}\lgr
    \Omega^{{\rm GL}_r,N}$ inducing a morphism  $\overline\Omega^{P,N}\lgr
    \overline\Omega^{{\rm GL}_r,N}$ which makes the diagram
  $$\def\mapdown#1{\Big\downarrow
  \rlap{$\vcenter{\hbox{$\scriptstyle#1$}}$}}
\matrix{
   P_{\!\!ad}^{N+1}/P_{\!\!ad}&\lgr&{\rm PGL}_r^{N+1}/{\rm PGL}_r&\cr
  \mapdown{}&&\mapdown{}\cr
&&\cr
     \overline\Omega^{P,N}&\lgr&
    \overline\Omega^{{\rm GL}_r,N}&\cr
  }$$
commutative. 
$\square$ 
\vskip2mm
Let $\pi_i$ and $\delta_i$ be respectively the $2(N+1)$-uples
$$((0,0),\cdots,(0,0),(1,0),(n-1,0),(0,0),\cdots,(0,0))$$ and 
$$((0,0),\cdots,(0,0),(0,0),(n,0),(0,0),\cdots,(0,0))$$ (where
$n-1$ and
$n$ are in the $(i+1)$-th position, $i\in{\Bbb Z}/(N+1){\Bbb Z}$), and
${\cal T}^{P,N}$ be the subtorus of $\smash{\vbox{\hbox{${\Bbb G}_m^{S^{P,N}}$}}}$ defined by the
$(N+1)$ equations
$\lambda_{\pi_i}=\lambda_{\delta_i}\quad (i\in {\Bbb Z}/(N+1){\Bbb Z})$.

The proof of the following proposition is now an easy calculation using the above facts~;
it will be left to the reader.
 \th P{ROPOSITION} 3.3
\enonce

 a) The morphism $$\displaylines{
\qquad P^{N+1}\times {\Bbb G}_m^{S^{P,N}}/(P\times {\Bbb
G}_m^{N+1})\longrightarrow
P^{N+1}/ P\times{\Bbb G}_m\hfill\cr
\hfill=\!= 
\{(g_0,\cdots,g_N,t)\in
P^{N+1}\times{\Bbb G}_m\mid g_0\cdots g_N=t\}\qquad}$$
defined by
$$(g_0,\cdots g_N,(\lambda_j)_j)\mapsto(
\lambda_{\pi_0}/\lambda_{\delta_0}g_0g_1^{-1},\cdots,\lambda_{\pi_N}/
\lambda_{\delta_N}g_Ng_0^{-1},\prod_i\lambda_{\pi_i}/\lambda_{\delta_i})$$
is an isomorphism

b) The morphism of tori ${\Bbb G}_m^{S^{P,N}} \longrightarrow{\Bbb G}_m^{S^{P,2}}$ associated
with an injection $\{j,j+1\}\hookrightarrow \{0,\cdots,N\}$ sends ${\cal T}^{P,N}$ to ${\cal
T}^{P,2}$

c) The composite $\Pi_i$ of the obvious morphism
$$\displaylines{\qquad\Omega^{P,N}\longrightarrow{\Bbb G}_m\backslash \prod_{i\in
S^{P,N}}[(\Lambda^{\vert i\vert}E/\sum_{j>i}\Lambda^jE)-\{0\}]\longrightarrow\hfill\cr
\hfill{\Bbb G}_m\backslash\Big[[(\Lambda^{\vert
\delta_i\vert}E/\sum_{j>\delta_i}\Lambda^jE)-\{0\}]\times[(\Lambda^{\vert
\pi_i\vert}E/\sum_{j>\pi_i}\Lambda^jE)-\{0\}]\Big]\qquad}$$
with the identifications
$$\displaylines{\qquad{\Bbb G}_m\backslash\Big[[(\Lambda^{\vert
\delta_i\vert}E/\sum_{j>\delta_i}\Lambda^jE)-\{0\}]\times[(\Lambda^{\vert
\pi_i\vert}E/\sum_{j>\pi_i}\Lambda^jE)-\{0\}]\Big]\hfill\cr
\hfill=\!={\Bbb
G}_m\backslash\Big[[(\Lambda^{\vert
\delta_i\vert}E)-\{0\}]\times[(\Lambda^{\vert
\pi_i\vert}E-\{0\}]\Big]\qquad}$$
and
$$\displaylines{\qquad{\Bbb G}_m\backslash\Big[[(\Lambda^{\vert
\delta_i\vert}E)-\{0\}]\times[(\Lambda^{\vert
\pi_i\vert}E-\{0\}]\Big]\hfill\cr
\hfill=\!=E_i\otimes[\Lambda^{n-1}E_{i+1}\otimes
(\Lambda^nE_{i+1})^{-1}]-\{0\}=\!=
{\goth gl}_n-\{0\}\qquad}$$
is ${\cal T}^{P,N}$-invariant and factors through $\goth P_{\rm rk\geq 1}$ . Its restriction 
to $P^{N+1}\times {\Bbb G}_m^{S^{P,N}}/( P\times {\Bbb
G}_m^{N+1})$ is just $$(g_0,\cdots g_N,(\lambda_j)_j)\mapsto(\lambda_{\pi_i}/\lambda_{\delta_i})
\,g_ig_{i+1}^{-1}$$
The ideal sheaves spanned by the minors of size $j$ of the matrix $\Pi_i$ are
invertible, $\forall\, 1\leq j\leq n-1$. So are the ideal sheaves spanned by the minors of size
$j$ of the matrix $A_i$ deduced from $\Pi_i$ via the obvious projection ${\goth P}\lgr{\goth
gl}_r$, $\forall\, 1\leq j\leq r-1$. In fact, this last matrix is a multiple of the one
obtained by composing the morphism $\Omega^{P,N}\lgr\Omega^{{\rm GL}_r,N}$ with the morphism
$\Pi_i^{\rm GL_r}$.

d) The kernel of the morphism ${\Bbb G}_m^{S^{P,N}}/\,{\Bbb G}_m^{N+1}\longrightarrow
 {\Bbb G}_m$
defined by
$(\lambda_j)_j\mapsto\prod_i\lambda_{\pi_i}/\lambda_{\delta_i}$ is a torus. This morphism 
extends to a
${\cal T}^{P,N}$ -invariant morphism 
$$t\ :\ {\cal A}^{P,N}\lgr{\Bbb A}^1\,.$$

e) The morphism $$(\Pi_0,\cdots,\Pi_N,t)\ :\ {\cal T}^{P,N}\backslash\Omega^{P,N}\longrightarrow
\goth P_{\rm rk\geq 1}^{N+1}\times{\Bbb A}^1_t$$ thus defined factors through $\mu^{P,N}$.

f) The morphism ${\cal T}^{P,N}\backslash\Omega^{P,N}\longrightarrow \mu^{P,N}$ is projective
and birational. It factors through $\mu_{\rm bl.}^{P,N}$ (cf. {\rm
2.3.6}\/).\quad$\square$\endth
\noindent
This achieves the construction of $\widetilde\mu^{P,N}={\cal T}^{P,N}\backslash\Omega^{P,N}$.
\vskip7mm
{\bf 4\pointir Complements and remarks}
\vskip 3mm 
4.1\pointir Faltings's original reduction principle.

Using the previous work  of L. Lafforgue on the compactification of ${\rm
PGL}_r^{N+1}/{\rm PGL}_r$ ([L~1]), it is also possible to obtain ${\rm
GL}_r^{N+1}$-equivariant log-smooth resolutions of
the schemes
$\mu^{r,N}$ ({\it cf}. 2.1) introduced by G. Faltings in [F 2].

This construction can  be  combined with Faltings's original  reduction principle [F
2] and one also obtains  log-smooth resolutions
$$\widetilde{\Bbb M}'_t\lgr{\Bbb M}_t$$
of the local models. In fact, the resolution $\widetilde{\Bbb M}_t\lgr{\Bbb M}_t$ constructed
in the present work (sections 2-3) factors through $\widetilde{\Bbb M}'_t\lgr{\Bbb M}_t$ ;
this reflects the fact that the morphism 
$$P_{\!\!ad}^{N+1}/P_{\!\!ad}\lgr {\rm PGL}_r^{N+1}/{\rm PGL}_r$$
extends to a morphism
$$\overline\Omega^{P,N}\lgr\overline\Omega^{{\rm GL}_r,N}$$
of Lafforgue's compactifications ({\it cf.} [L 2]).

The reason why we developed our variant is that it seems to be a little bit more symmetric
({\it cf.} remark 4.2.3).

4.1.1. {\it Remark}: in [F 2], Faltings constructs a proper birational morphism 
$$Y\lgr X=\mu^{r,N}$$
({\it cf}. [F 2], p. 25) and checks that for $r\leq 3$ this morphism is a  
${\rm GL}_r^{N+1}$-equivariant log-smooth resolution of $\mu^{r,N}$. Combined with his
reduction principle, this result yields resolutions of the local models for $r\leq 3$ in the
unitary case and $g\leq 3$ in the symplectic one. I do not know the relationship between
Faltings's resolutions and the ones we obtain by using Lafforgue's work. I also do not know the
relationship of all these resolutions  with the one constructed by de Jong for $g=2$ ([dJ
1]) and with the one obtained in [G] for $g=3$.
\vskip2mm
4.2\pointir
In some cases the local model ${\Bbb M}_t$ is in fact endowed with an action of a
non-connected group-scheme ${\cal K}^{\rm ext.}$  with neutral component ${\cal K}$. More
precisely, in the unitary case, this occurs when a cyclic permutation of $\{1,\cdots, N+1\}$
fixes the partition $(d_i)_{1\leq i\leq N+1}$ of $n$ and in the symplectic case, this occurs
when the permutation $(N,N-1,\cdots,1)$ fixes the partition $(d_i)_{1\leq i\leq N}$ of $g$. We
shall see that in these cases the resolution \smash{$\widetilde {\Bbb M}_t$} is also endowed with
an action of
${\cal K}^{\rm ext.}$.
\vskip2mm
4.2.1\pointir The unitary case.
Let $\sigma\ :\ i\lgr i+s\quad(s\in {\Bbb Z}/(N+1){\Bbb Z})$ be a cyclic permutation of
$\{1,\cdots,N+1\}={\Bbb Z}/(N+1){\Bbb Z}$ fixing $(d_i)_i$ and $S$ be the order of $\sigma$.
The group ${\Bbb Z}/S{\Bbb Z}=\,<\sigma>$ acts on ${\cal K}$, on  ${\Bbb M}_t({\rm
GL}_n\,,\,r,(d_i)_i)$ and on $\mu^{P,N}$ through the permutation of  the factors $V_i$'s induced
by
$\sigma$. This defines an action of ${\cal K}^{\rm ext.}:={\cal K}\rtimes ({\Bbb Z}/S{\Bbb
Z})$ on
${\Bbb M}_t({\rm GL}_n\,,\,r,(d_i)_i)$. It is easily seen that  the resolution
$\widetilde\mu^{P,N}\lgr \mu^{P,N}$   is ${\Bbb Z}/S{\Bbb Z}$-equivariant:
in fact Lafforgue's compactifications are even endowed with an action of the symmetric group
${\goth G}_{N+1}$ ({\it cf.} [L 1,2], 1.b). Hence the resolution
$\widetilde{\Bbb M}_t({\rm GL}_n\,,\,r,(d_i)_i)\lgr{\Bbb M}_t({\rm
GL}_n\,,\,r,(d_i)_i)$ induced by $\widetilde\mu^{P,N}$ is ${\cal K}^{\rm ext.}$-equivariant.
\vskip2mm
4.2.2\pointir The symplectic  case 

Let $\sigma$ be the permutation $(N,N-1,\cdots,1)$ of $\{1,\cdots,N\}$. We also denote by the
same letter the permutation $(N,N-1,\cdots,1,N+1)$ of $\{1,\cdots,N+1\}$. The group
${\Bbb Z}/2{\Bbb Z}=\,<\sigma>$ acts on ${\cal K}$, on ${\Bbb M}_t\,({\rm
Sp}_{2g}\,,(d_i)_i)$ and on $\mu^{P,N}$, respectively by letting the generator $\sigma$ act via
$$\leqalignno{(\gamma_i)_i&\mapsto (J^{-1}\ {}^t\gamma_{\sigma'(i)}^{-1}J)_i\,,\cr
(\omega_i)_i&\mapsto(\omega_{\sigma'(i)}^\perp)_i\cr
&&{\rm and}\hfill\cr
(\Pi_i)_i&\mapsto (J^{-1}\ {}^t\Pi_{\sigma(i)}J)_i }$$
(where $J$ is the matrix (2.3),  the chain
$$\matrix{V_N^\vee{} &\buildrel{{}^t\alpha_0}\over\lgr
&V_0^\vee&\buildrel{{}^t\alpha_1}\over\lgr& V_1^\vee&\buildrel{{}^t\alpha_2}\over\lgr&
\cdots&\buildrel{{}^t\alpha_N}\over\lgr &V_N^\vee\cr
\hphantom{V_0}&\hphantom{\buildrel{\alpha_0}\over\lgr} &\hphantom{V_N}
&\hphantom{\buildrel{\alpha_N}\over\lgr}&
\hphantom{V_{N-1}}&\hphantom{\buildrel{\alpha_{N-1}}\over\lgr}&\hphantom{\cdots}&
\hphantom{\buildrel{\alpha_1}\over\lgr}
&\hphantom{V_0}}$$\vskip-5mm\noindent   is identified to the chain
$$\matrix{V_0&\buildrel{\alpha_0}\over\lgr &V_N &\buildrel{\alpha_N}\over\lgr&
V_{N-1}&\buildrel{\alpha_{N-1}}\over\lgr&\cdots&\buildrel{\alpha_1}\over\lgr &V_0\cr
} $$
\eject
by the perfect pairing on ${\Bbb Z}[t,t^{-1}]^{2g}$ defined by the matrix
$$\pmatrix{0&{\rm Id}_g\cr
t.{\rm Id}_g&0}J$$
and $\sigma'$ is the permutation $(N,N-1,\cdots,0)$ of
$\{0,\cdots,N\}$). This defines an action of ${\cal K}^{\rm ext.}:={\cal K}\rtimes ({\Bbb
Z}/2{\Bbb Z})$ on
${\Bbb M}_t({\rm Sp}_{2g}\,,\,(d_i)_i)$. One can see that if the resolution
$\widetilde\mu^{P,N}\lgr \mu^{P,N}$   is ${\Bbb Z}/2{\Bbb Z}$-equivariant (and this will be
the case for the resolutions obtained in the third section from Lafforgue's
compactifications: the automorphism $g\mapsto J^{-1}{}^tg^{-1}J$ of $P$ induces an
automorphism of Lafforgue's compactifications, {\it cf.} [L 2], 1.b), the resolution
$\widetilde{\Bbb M}_t({\rm Sp}_{2g}\,,\,(d_i)_i)\lgr{\Bbb M}_t({\rm
Sp}_{2g}\,,\,(d_i)_i)$ induced by $\widetilde\mu^{P,N}$ is ${\cal K}^{\rm ext.}$-equivariant.
\vskip 2mm
4.2.3. {\it Remarks}:

1) The reason why we have developed our variant and used Lafforgue's
generalization [L 2] of his previous work is that in the symplectic case, the symmetry $\sigma$
of ${\Bbb M}_t({\rm
Sp}_{2g}\,,\,(d_i)_i)$ does not (or at least not obviously) lift to an automorphism of
\smash{$\widetilde{\Bbb M}'_t({\rm Sp}_{2g}\,,\,(d_i)_i)$}.

2) The   resolution of ${\Bbb M}_p({\rm
Sp}_{4}\,,\,(1,1))$ constructed by de Jong in [dJ 1] is also 
${\cal K}_p^{\rm ext.}$-equivariant.  I do not know wether this is the
case for the resolution of
${\Bbb M}_p({\rm Sp}_{6}\,,\,(1,1,1))$ constructed in [G].

\vskip 12mm
{\bf Bibliography}
\vskip 5mm
\leftskip=\parindent
\parindent=0mm
\parskip=2mm

\leavevmode\llap{[CN]} C. N. Chai and P. Norman\pointir Bad reduction of the Siegel moduli
scheme of genus two with $\Gamma_0(p)$-level structure, {\sl Amer. Jour of Math.} {\bf 112}
(1990) pp. 1003--1071

\leavevmode\llap{[dJ 1]} A. J. de Jong\pointir Talk in Oberwolfach (July 1992).

\leavevmode\llap{[dJ 2]} A. J. de Jong\pointir The moduli spaces of principally
polarized abelian varieties with $\Gamma_0(p)$-level structures,
{\sl J. Alg. Geometry} {\bf 2} (1993) pp. 667--688

\leavevmode\llap{[DP]} P. Deligne and G. Pappas\pointir Singularit\'es des espaces de modules
de Hilbert en les caract\'eristiques divisant le discriminant, {\sl Comp. Math} {\bf 90} (1994)
pp. 59--79

\leavevmode\llap{[DR]} P. Deligne and M. Rapoport\pointir Modules des courbes elliptiques, {\sl
Modular functions of one variable} {\bf SLN 349}, Springer Verlag (1973) 

\leavevmode\llap{[F 1]} G. Faltings\pointir Explicit resolution of singularities of moduli
spaces, {\sl Journ. reine angew. Math.} {\bf 483} (1997) pp.183--196

\leavevmode\llap{[F 2]} G. Faltings\pointir Toro\"\i dal resolutions for some matrix
singularities, {\sl preprint} (1999)

\leavevmode\llap{[G]} A. Genestier\pointir Un mod\`ele semi-stable de la vari\'et\'e de Siegel de
genre 3 avec structures de niveau de type $\Gamma_0(p)$,
{\sl preprint (Universit\'e Paris-Sud, 98--72)}, to appear in {\sl Compositio Mathematica}

\leavevmode\llap{[KKMS]} G. Kempf, F. Knudson, D. Mumford and B. Saint-Donnat\pointir
{\sl Toro\"\i dal embeddings~I,} {\bf SLN 339}, Springer Verlag (1971)

\leavevmode\llap{[L 1]} L. Lafforgue\pointir Pavages des simplexes, sch\'emas de graphes
recoll\'es et compactification des ${\rm PGL}_r^{n+1}/{\rm PGL}_r$,
{\sl Inventiones mathematicae} {\bf136} (1999) pp. 233--271

\leavevmode\llap{[L 2]} L. Lafforgue\pointir Compactification des ${\rm PGL}_r^{n+1}/{\rm PGL}_r$
et restriction des scalaires \`a la Weil, {\sl preprint (Universit\'e Paris-Sud, 98--47)}, to
appear in {\sl Annales de l'institut Fourier}

\leavevmode\llap{[R]} M. Rapoport\pointir On the bad reduction of Shimura varieties, 
{\sl Automorphic forms, Shimura varieties and $L$-functions vol. II (\'ed. by L. Clozel 
and J. S. Milne),
Persp. in Math. } {\bf 11} (1990) pp. 253--321

\leavevmode\llap{[RZ]} M. Rapoport and Th. Zink\pointir Period Spaces for $p$-divisible
Groups, {\sl Annals of Math. Studies} {\bf 141} (Princeton University Press, 1996)

\vskip 5mm
\leavevmode\hfill\vbox{\hbox{Alain Genestier}
\hbox{UMR 8628, CNRS}
\hbox{}
\hbox{Universit\'e Paris-Sud}
\hbox{Math\'ematique, B\^at. 425}
\hbox{91405 Orsay (France)}
\hbox{}
\hbox{e --mail:}
\hbox{\tt Alain.Genestier@math.u-psud.fr}}
\bye